\documentclass[review]{elsarticle}

\usepackage{hyperref}
\usepackage{amsmath,amsthm}
\usepackage{amssymb}

\begin{document}

\begin{frontmatter}

\title{Reflecting stochastic dynamics of active-passive populations with applications in operations research and neuroscience}


\author[mymainaddress]{Thi Kim Thoa Thieu \corref{mycorrespondingauthor}}
\cortext[mycorrespondingauthor]{Corresponding author}

\author[mymainaddress,mysecondaryaddress]{Roderick Melnik}


\address[mymainaddress]{M3AI Laboratory, MS2Discovery Interdisciplinary Research Institute, Wilfrid Laurier University, \\75 University Ave W, Waterloo, Ontario, Canada N2L 3C5}
\address[mysecondaryaddress]{BCAM - Basque Center for Applied Mathematics, Bilbao, Spain \\Email: \{tthieu, rmelnik\}@wlu.ca}

\begin{abstract}
Stochastic dynamic models have been extensively used for the description of processes with uncertainties arising in the operations research, behavioral sciences, and many other application areas.  A large class of the problems from these domains is characterized by the necessity to deal with several distinct groups of populations, which are usually labeled as \textquotedblleft active\textquotedblright  and \textquotedblleft passive\textquotedblright. Motivated by important applications of queueing networks and neuroscience, the main focus of the present work is on the analysis of reflecting stochastic dynamics of such mixed populations. We develop a general mathematical modeling framework to describe the reflecting stochastic dynamics for active-passive populations. The analysis of this model is carried out via a combination of low- and high-delity results obtained from the solution of the underlying coupled system of SDEs
and from the simulations with a statistical-mechanics-based lattice gas model, where we employ a kinetic Monte Carlo procedure. We provide details of the queueing theory and neuronal models and discuss a relationship between reflecting SDEs and a model of queueing theory via a limit theorem. Furthermore, we present several representative numerical examples, and discuss an intrinsic interconnection between active and passive particles in the underlying stochastic process. Finally, possible extensions of the proposed methodology have been highlighted.
\end{abstract}


\end{frontmatter}


\section{Introduction}
\label{Intro}

Complex behavioral systems are ubiquitous in many areas of science, including biology, neuroscience, engineering, as well as in social sciences and operations research in general.
The development of more realistic complex behavioral models offers many challenging questions to science and technology.
Hence, in the modelling of such mathematical models, the uncertainty quantification plays an important role for the setting of mathematical formulation of the problems \cite{Bellomo2017}. However, in terms of the mathematical formulation of the corresponding problems, the overall dynamics of such systems should, in most cases, be considered in confined domains, possibly with obstacles. For this particular class of stochastic dynamic models, the boundary conditions must be specified and the frequent choices are absorbing and reflecting or even mixed boundary conditions. Therefore, we turn our attention to the recent developments provided in \cite{Thieu2020-3,Zhang2016,Battiston2020, Thieu2021} that bring us to the study of a system of reflecting stochastic dynamics, namely Skorokhod-type stochastic differential equations (SDEs), modelling the dynamics of active-passive populations.
Furthermore, operations research that includes queueing theory models is one of important applications of such reflecting stochastic dynamics. The class of queueing theory models is important in many applied fields, including the emergency medical aid systems, passenger services
in air terminals, bio-social systems and neuroscience, as well as various types of other systems and processes. The models of queueing theory and those based on systems of SDEs are frequently studied separately, although the links between them have been explored by a number of authors. In \cite{Ward2003}, the authors considered a single-server queue with a Poisson arrival process and exponential processing times in which each customer independently reneges after an exponentially distributed amount of time.  The authors in \cite{Ward2005} discussed a single-server queue with a renewal arrival process and generally distributed processing
times in which each customer independently reneges if service has not begun within a distributed amount of time. The topic of heavy traffic limit approximations in queueing network models has been investigated in \cite{Williams1998, Ramanan2003, Kruk2011}. In \cite{Ramanan2006}, the author considered the extended Skorokhod problem (ESP) and associated extended Skorokhod map (ESM) that enabled a pathwise construction of reflected diffusions. In \cite{Thieu2020-3}, the authors proved the well-posedness of a coupled system of Skorohod-like SDEs modeling the dynamics of active-passive pedestrians. A discretization scheme for reflected SPDEs driven by space-time white noise through systems of reflecting SDEs has been discussed in \cite{Zhang2016}. We note also \cite{Pilipenko2014} where a discussion and references on a relationship between reflecting SDEs and some models of queueing theory were provided. 

 On the other hand, our better understanding of the activities of neuronal cells in nervous systems is one of the major current challenges.
 The neural activity is intrinsically noisy and the specific neural activity in living organisms can be treated via the paradigm of complex systems, so that stochastic dynamics with various types of uncertainties have to be incorporated into the models (see e.g. in \cite{Laing2010}). Furthermore, since the biophysical dynamics of neurons (e.g. membrane potential and synaptic processes) are complicated and hard to measure, then stochastic tools can be useful in the study of such dynamics of neurons. Up to date, a number of relevant results are available in the topic of stochastic dynamics of neurons. One of the
most prominent examples for this class of models is given by the Hodgkin-Huxley
model \cite{Hodgkin1952} that has been introduced in 1952 to describe the membrane potential in the
squid giant axon. A two-compartment neuronal model has been reported in \cite{Lansky1999}.   
In \cite{Reutimann2003}, the authors proposed event-based strategy for efficiently simulating large networks of simple
model neurons. In \cite{Tusbo2012}, the authors have shown an
alternative hypothesis belonging to the maximization of mutual information for neural code from the neuronal activities recorded juxtacellularly in the sensorimotor cortex of
behaving rats. An inverse first passage time method for a two dimensional Ornstein-Uhlenbeck process and its applications in neuroscience have been discussed in \cite{Civallero2019}. Finally, the authors in \cite{Hsu2021} have provided a review of the  current challenges and further progress in building stochastic models for single-cell data. 

Motivated by the study of a system of Skorokhod-type SDEs together with applications in queueing theory models and neuroscience, our representative examples here are pertinent to (i) the modelling of customer services in a queueing theory setting in general and (ii) neural dynamics in a cell's membrane potential model. 


In this paper, we introduce two reflecting stochastic dynamic models for active-passive populations. In particular, in the first model, we study the stochastic dynamics of active-passive populations that can be represented as customer services in a queueing theory setting. Specifically, the system consists of a crowd of customers based on two distinct populations: (a) active customers who have priority and valid accesses to the servers and leave the domain immediately after reaching the servers, and (b) passive customers who do not have priority and valid accesses to the servers. In what follows we show that the model of queueing theory converges to a system of reflected SDEs via a limit theorem. In the second model, we study a cell's membrane potential model of two distinct populations of ions: $\text{Na}+$ and $\text{Cl}-$. The main question we ask, while applying this model, is \textquotedblleft How to capture the behaviour of many ions when the action potential reaches peak and then goes down to the hyperpolization state?\textquotedblright. To address this question we consider a system of SDEs of Skorokhod type with reflecting boundary conditions based on similar ideas of a particular case of our queueing theory model to capture the behaviour of two different ions in a cell's membrane potential model. The analysis of this model is carried out by combining low- and high-fidelity results obtained from the solution of the underlying coupled system of SDEs and from the simulations with a statistical-mechanics-based lattice gas model, where we employ a kinetic Monte Carlo
procedure.

The rest of this paper is organized as follows. In Section \ref{model}, we provide details of the model and consider an application of reflecting stochastic dynamics in queueing theory. A system analysis based on a multi-fidelity approach has also been presented in this section. In Section \ref{neuro-model}, we focus on an application of reflecting stochastic dynamics in neuroscience. The hyperpolarization phenomenon and its influence on the dynamics of neuronal cells have been discussed in detail. Finally, concluding remarks are provided in Section \ref{concluding-remarks}.

\section{An application in queueing theory}\label{model}
In this section, we consider a system of reflecting stochastic dynamics for an active-passive customers service. In particular, we define the model, develop its implementation, and provide numerical examples.
\subsection{Model description}\label{model-1}
We consider an arriving customer service for active and passive customers. The geometry is a square lattice $\Lambda:=\{1,\dots,L\}\times\{1,\dots,L\}\subset \mathbb{Z}^2$ of side 
length $L$, where $L$ is an odd positive integer number. $\Lambda$ will be referred in this context as {\em room}.  
An element $x=(x_1,x_2)$ of the room $\Lambda$ is called \emph{site} 
or \emph{cell}.
Two sites $x,y\in\Lambda$ are said \emph{nearest neighbors} if and only 
if $|x-y|=1$. In addition, $N_A$ is the total number of active customers, $N_P$ is the total number of passive customers with $N:=N_A+N_P$ and $N_A, N_P, N \in \mathbb{N}$. 

We assume that an exit door located on the top row of the geometry, and that there are $\omega$ servers at the exit door. Every active customer's target is to reach these servers and to leave the geometry immediately after that, while the passive customers cannot leave the domain due to invalid access to the servers. 
\begin{figure}[h!]
	\centering
	\includegraphics[width=0.8\textwidth]{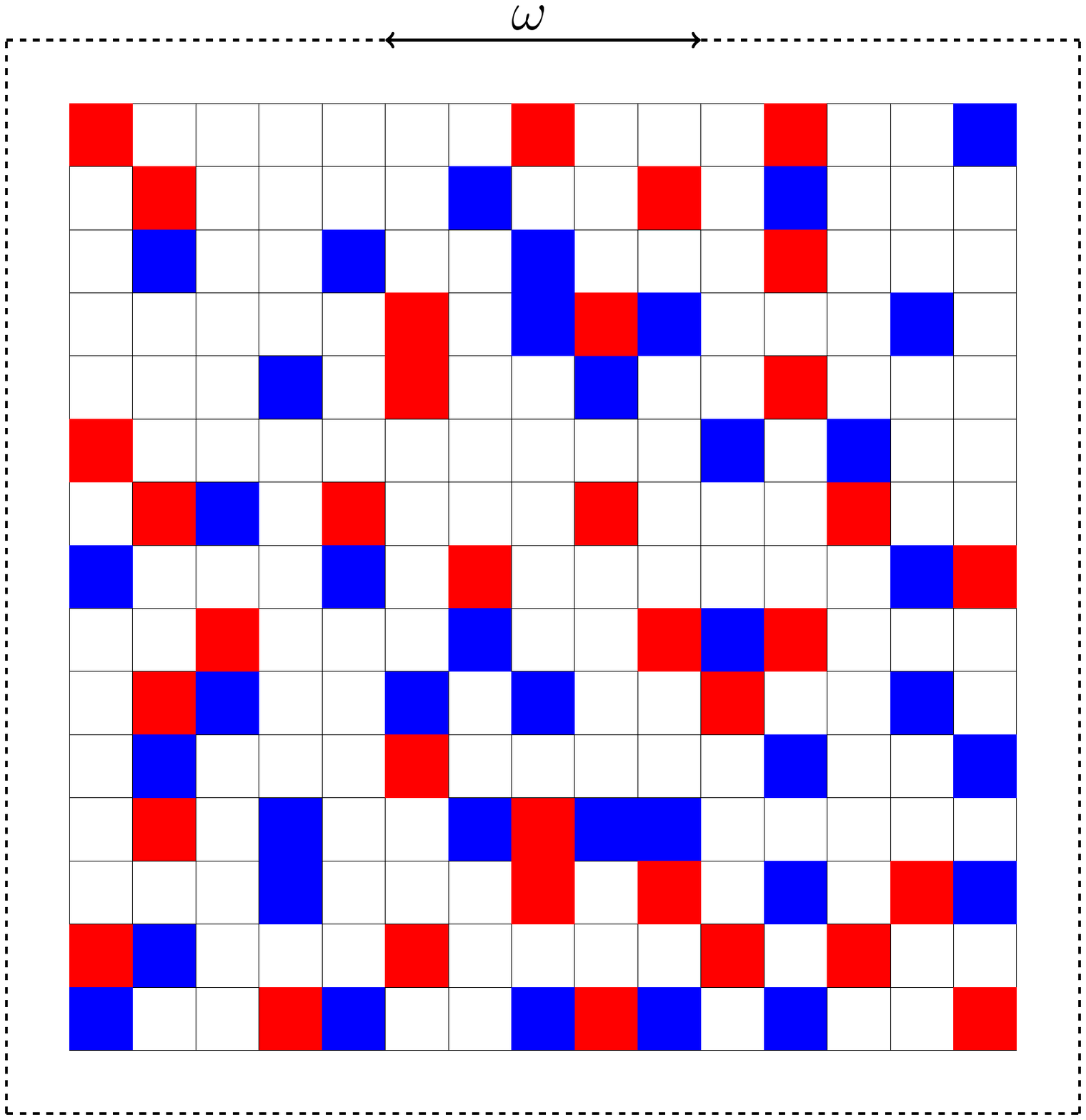}
	\vspace*{-40mm}
	\caption{(Color online) Schematic representation of our lattice model. Blue and red 
		squares denote passive and active customers, respectively, while the white squares within the geometry represents the empty spots. The thick dashed line surrounding a large
		fraction of the grid denotes the presence of reflecting boundary conditions. The exit door is located in presence of the arrow, with its width equal to $\omega$.}
	\label{fig:fig0}
\end{figure}
This behavioral system in the context of queueing theory can be viewed as an $M/M/\omega/N$\footnote[1]{The notation $M/M/\omega/N$ means that the queueing model contains $\omega$ serves and $N$ places for waiting, while $M$ represents Poisson arrival process (i.e., exponential inter-arrival times) see e.g. \cite{Kalashnikov1994}.} queueing system with the following assumptions:
\begin{itemize}
	\item [a)] there are $N$ customers that request services,
	\item [b)] the maximum number of requests in a system equals $\omega$,
	\item [c)] the arrival times of active and passive customers are exponential i.i.d random variables with intensities $\alpha$ and $\mu$, respectively.
\end{itemize}

Next, we highlight key steps in our model construction.
\begin{itemize}
	\item[i.] The active population of the system is described as the following process:
	
	\begin{align}\label{active-eqn}
	X_{A_i}^\alpha(t) = X_{A_i}^\alpha(0) + Y_{A_i}^{\alpha}(t) + \Phi_{A_i}^{\alpha}(t),
	\end{align}
	where $X_{A_i}^\alpha(t)$ denotes the position of the customer $i \in \{1,\ldots, N_A\}$ at time $t\geq 0$. On the other hand, $Y_{A_i}^{\alpha}$
	can be interpreted as the arrival times of active customers, $\Phi_{A_i}^{\alpha}(t)$ is the cumulative lost service capacity over $[0,t]$.
	\item[ii.] The passive population of the system is described as the following process:
	
	\begin{align}\label{passive-eqn}
	X_{P_j}^\beta(t) = X_{P_j}^\beta(0) + Y_{P_j}^{\beta}(t) + \Phi_{P_j}^{\beta}(t),
	\end{align}
	where $X_{P_j}^\beta(t)$ denotes the position of the customer $j \in \{1,\ldots, N_P\}$ at time $t \geq 0$. Similarly, $Y_{P_j}^{\beta}$
	can be interpreted as the arrival times of passive customers, $\Phi_{P_j}^{\beta}(t)$ is the cumulative lost service capacity over $[0,t]$. 
\end{itemize}

Motivating by \cite{Pilipenko2014} (cf. Section 3.3), we consider a limit theorem for a model of queueing theory that leads us to a system of
reflecting SDEs.
Using the central limit theorem, we have

\begin{align}\label{limit-1}
\frac{Y_{A_i}^{\alpha}(t) - \alpha t}{\sqrt{\alpha}} \to w_{A_i}(t) \text{ when } \alpha \to \infty,
\end{align}
and
\begin{align}\label{limit-2}
\frac{Y_{P_i}^{\beta}(t) - \beta t}{\sqrt{\beta}} \to w_{P_j}(t) \text{ when } \beta \to \infty
\end{align}
in $\Lambda$.
Similarly, we also have

\begin{align}\label{Skorohod-1}
\left(\frac{X_{A_i}^\alpha(t)}{\sqrt{\alpha}},\frac{\Phi_{A_i}^\alpha(t)}{\sqrt{\alpha}} \right) \to (\xi_{A_i}(t), \phi_{A_i}(t)) \text{ when } \alpha \to \infty,
\end{align}
and
\begin{align}\label{Skorohod-2}
\left(\frac{X_{P_j}^\beta(t)}{\sqrt{\beta}},\frac{\Phi_{P_j}^\beta(t)}{\sqrt{\beta}} \right) \to (\xi_{P_j}(t), \phi_{P_j}(t)) \text{ when } \beta \to \infty
\end{align}
in $\Lambda$.

As a result, the pairs $(\xi_{A_i}(t), \phi_{A_i}(t))$ and $(\xi_{P_j}(t), \phi_{P_j}(t))$ are the solutions of the following Skorokhod type equations:

\begin{align}
\xi_{A_i}(t) &=  w_{A_i}(t) + \phi_{A_i}(t), \\
\xi_{P_j}(t) &=  w_{P_j}(t) + \phi_{P_j}(t).
\end{align} 

In the lattice-based numerical implementation of this model, each customer is represented as a particle. The resulting coupled particle system is interacting via the site exclusion principle, i.e. each site of the lattice can be occupied only by a single particle.


\subsection{Numerical results}\label{KMC}
Our representative examples are reported here for two different species of
\emph{active} and \emph{passive} customers, moving inside $\Lambda$
(we use in the notation the symbols A and P to 
respectively refer to them).
Note that the sites 
of the external boundary of the room, i.e. the sites  
$x\in\mathbb{Z}^2\setminus\Lambda$, is such that there exists 
$y\in\Lambda$ nearest neighbor of $x$ which
cannot be accessed by the customers. 
We call the state of the system
\emph{configuration} $\eta\in\Omega=\{-1,0,1\}^\Lambda$ 
and 
we shall say that the site $x$ is 
\emph{empty} if $\eta_x=0$,
\emph{occupied by an active customer} if $\eta_x=1$,
and
\emph{occupied by a passive customer} if $\eta_x=-1$.
The number of active (respectively, passive) 
customers in the configuration $\eta$ 
is given by 
$n_{\text{A}}(\eta)=\sum_{x\in\Lambda}\delta_{1,\eta_x}$
(resp.\ $n_{\text{P}}(\eta)=\sum_{x\in\Lambda}\delta_{-1,\eta_x}$), 
where $\delta_{\cdot,\cdot}$ is Kronecker's symbol.
Their sum is the total number of particles in the configuration $\eta$. The overall dynamics of our system is incorporated in the continuous time Markov chain $\eta(t)$ on $\Omega$ with rates $c(\eta,\eta')$, for the detailed definitions of the rates $c(\eta(t),\eta)$, see, e.g., \cite{Cirillo2019, Cirillo2020}. 

We took the inspiration from the population dynamics in \cite{Cirillo2019, Colangeli2019} to model the dynamics of our active and passive customers. Using the descriptions based on a simple exclusion process, the dynamics in the room is modeled as follow: the passive customers perform a 
symmetric simple exclusion dynamics on the whole lattice, whereas active customers is subject to a drift, guiding particles towards the exit door.

The numerical results reported in this section are obtained by using the kinetic Monte Carlo (KMC) method. In particular,
we 
simulate the presented model by using the 
following scheme: we 
extract at time $t$ an exponential random time $\tau$ as a function with parameter which has 
the total rate
$\sum_{\zeta\in\Omega}c(\eta(t),\zeta)$, then set the time $t$ equal to $t+\tau$.
Next, we select a configuration using the probability 
distribution 
$c(\eta(t),\eta)/\sum_{\zeta\in\Omega}c(\eta(t),\zeta)$
and set $\eta(t+\tau)=\eta$ (for the detailed definitions of the rates $c(\eta(t),\eta)$, see, e.g., \cite{Cirillo2019,Cirillo2020}).

This numerical scheme has initially been studied in \cite{Cirillo2019}, where the authors implemented a version of KMC to analyze the pedestrian escape from an obscure room by using a lattice gas model with two species of particles. Further validation of the numerical methodology used here was reported in \cite{Cirillo2020}, based on a counterflow pedestrian model. Note that the overall dynamics of our system is based on a continuous time Markov chain, i.e. the process will change state according to an exponential random variable and then move to a different state as specified by the probabilities of a stochastic matrix, together with a simple exclusion process. On the other hand, we use a statistical-mechanics-based lattice gas framework where we employ a kinetic Monte Carlo procedure to simulate the queueing theory model described in earlier in this section. In general, Monte Carlo statistical methods, particularly those based on Markov chains, provide a possible way to sample the distribution of the input random variable. Moreover, the reflecting stochastic queueing theory model we discussed here is a low-delity
modelling methodology for general behavioral systems. In order to approximate this
reflecting stochastic dynamics problem, we combine the coupled system of
SDEs with a high-delity modelling method by using a statistical-mechanics-based lattice
gas approach. Therefore, the methodology used in our model can be considered as a multi-fidelity approach to statistical inference (e.g. in \cite{Peherstorfer2018-2,Robert2004}).

Here, we consider the system defined in Section \ref{model-1} for $L=60$, $\omega=20$. All of the simulations are obtained starting with the same initial configuration chosen once for all by distributing the customers at random with uniform probability distribution. All other necessary details are provided in the figure captions.

\begin{figure}[h!]
	\centering
	\begin{tabular}{ll}
		\includegraphics[width=0.45\textwidth]{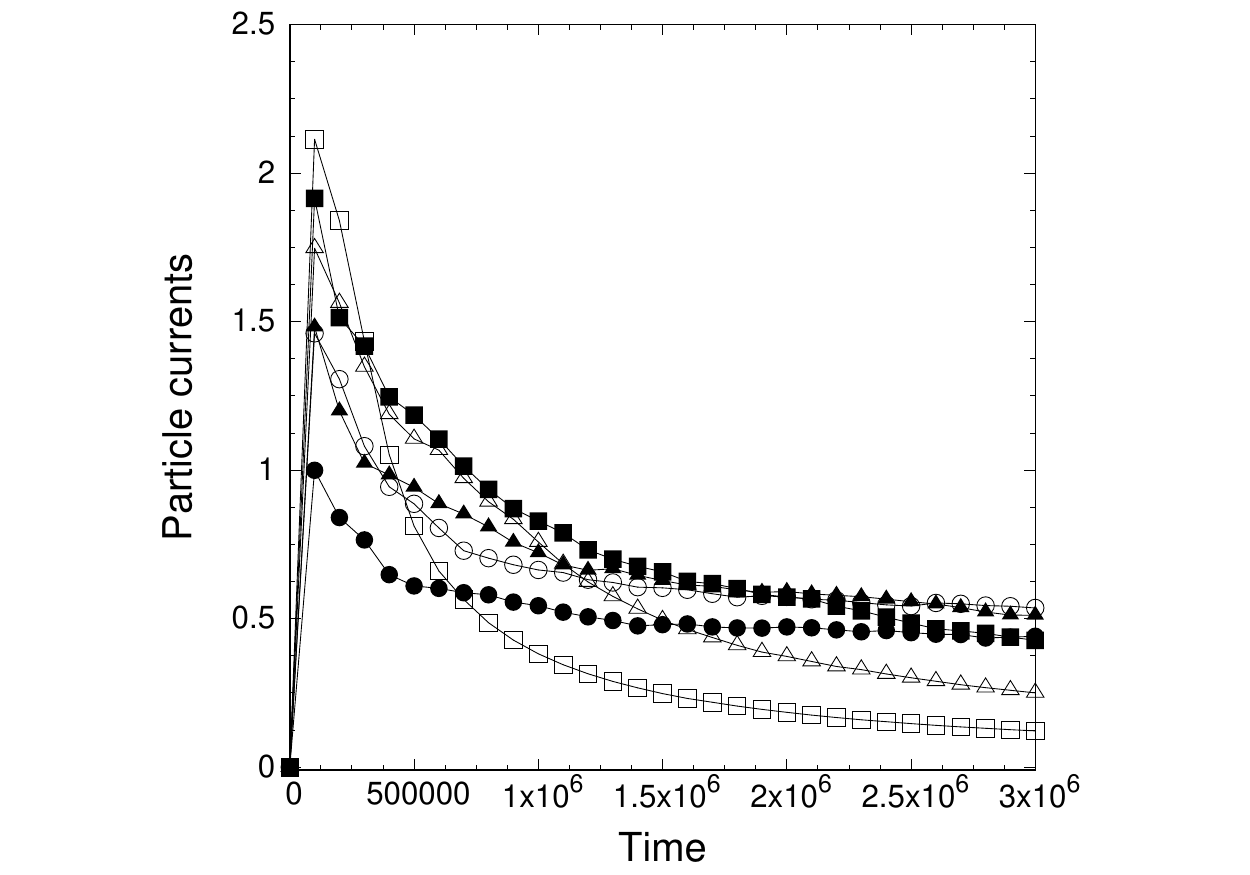}&
		\includegraphics[width=0.45\textwidth]{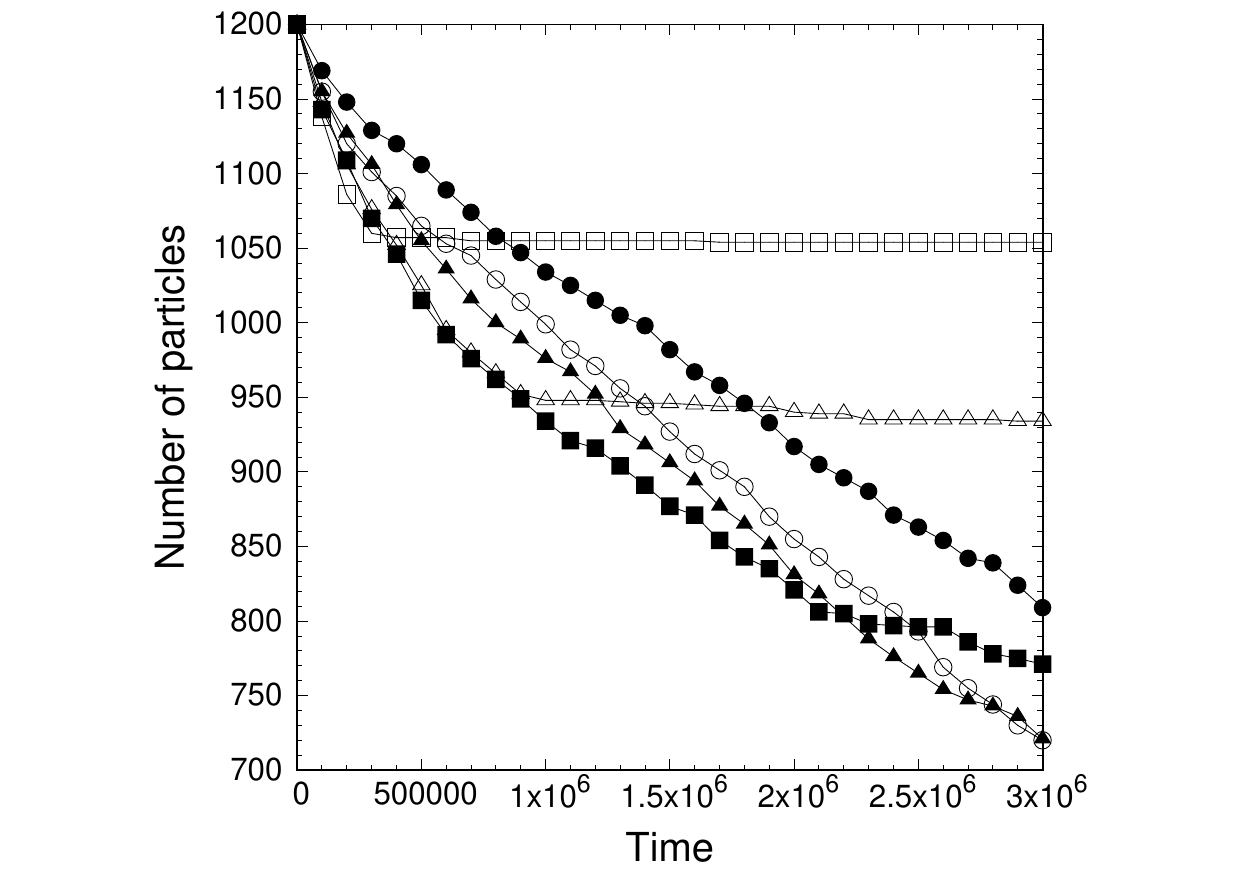}
	\end{tabular}
	\caption{Left panel: Evolution of the current as a function of time for active-passive customers and active customer exit times for $L=60$, $N_A=N_P=1200$ (empty symbols) and $N_A=1200, N_P=0$ (solid symbols) with $\varepsilon=0.1$ (circles), $\varepsilon=0.3$ (triangles), $\varepsilon=0.5$ (squares). Right panel: Behavioral pattern of the crowd for $L=60$, $N_A=N_P=1200$ (empty symbols) and $N_A=1200, N_P=0$ (solid symbols) with $\varepsilon=0.1$ (circles), $\varepsilon=0.3$ (triangles), $\varepsilon=0.5$ (squares).}
	\label{fig:fig1}
\end{figure}

The main representative numerical results of our analysis here are shown in Fig. \ref{fig:fig1}. We have plotted the total number of customers and the customer currents (outgoing fluxes) as functions of time (see e.g. in \cite{Colangeli2019}), for different values of $\varepsilon$ ($\varepsilon \geq 0$ is the drift quantity). Note that the current of active customers is defined in the infinite time limit by the
ratio between the total number of active customers, that in the interval $(0, t)$ passed
through the servers to leave the geometry, and the time $t$. Among other observations, our numerical results have demonstrated that the residence time of passive customers increases the residence time of active customers. Furthermore, we have observed that \textquotedblleft too smart\textquotedblright \ active customers increase their exit times. By the notion of \textquotedblleft too smart\textquotedblright, applied to active customers, we mean that when the values $\varepsilon$ are large enough, active customers move too quickly towards the exit door. In the left panel of Fig. \ref{fig:fig1}, the dynamics in the presence of passive customers increases the residence time of active customers. The current of active customers in the dynamics of only active customers is larger than in the case of both active and passive customers in the geometry for the case of $\varepsilon=0.1$, $\varepsilon=0.3$ and $\varepsilon=0.5$. To explain the observed effect, we examine the corresponding active customer exit times in the right panel of Fig. \ref{fig:fig1}. Due to the simple exclusion principle, the presence of passive customers can be seen as an obstacle in the geometry, explaining an increase in the waiting time of active customers. However, when increasing the value of $\varepsilon$ from $0.1$ to $0.3$ and $0.5$, the current of active customers in the dynamics of both active and passive customers with $\varepsilon=0.5$ is smaller than the current of the same dynamics with $\varepsilon=0.3$. Similarly, the current of active customers in the dynamics of both active and passive customers with $\varepsilon=0.3$ is smaller than the current of the same dynamics with $\varepsilon=0.1$. This effect is interesting. In general, when we increase the values of drift quantity, the current of active customers should be increasing. Looking at the corresponding active customer exit times in the right panel of Fig. \ref{fig:fig1}, we see the same effect where the residence time of active customers in the case of $\varepsilon=0.3$  is larger than in the case of $\varepsilon=0.1$. This is due to the fact that an increase in the value of $\varepsilon$ lets the active customers become \textquotedblleft too smart\textquotedblright. Hence, the active customers go quickly to the exit door and accumulate there together with the presence of passive customers in the geometry that makes the queue longer and increase the exit time of the active customers. This phenomenon could also be observed in the dynamics of only active customers where the current of active customers with $\varepsilon=0.5$ is larger than the current with $\varepsilon=0.3$ up to a certain value. After that it becomes smaller than in the case of $\varepsilon=0.3$. Similar effect can be observed in the corresponding active customer exit times in the right panel of Fig. \ref{fig:fig1}. Finally, it is worth noting that  \textquotedblleft too smart\textquotedblright \ active customers increase their residence times. As a closing note for this discussion, the presence of passive customers increase the time residence of active customers in the system. As further research, this numerical investigation would be interesting not only in the queueing theory models but also in the applications of neural networks, social networks and so on. As we mentioned above, our reflecting stochastic dynamics for active-passive populations can be applied also in a neuronal model, which we consider next. 

\section{An application in neuroscience}\label{neuro-model} 

In this section, we study a system of reflecting stochastic dynamics of a cell's membrane potential model for two distinct ions $\text{Na}+$ and $\text{Cl}-$. Specifically, we define the model and provide numerical examples.
\subsection{Model description}
For the functionality of a nervous system, neurons must be able to send and receive signals (see, e.g., \cite{Rye2016}). In particular, each neuron has a charged cellular membrane that causes a voltage difference between the inside and the outside. This charge of the membrane can change in response to neurotransmitter molecules released from other neurons and environmental stimuli. Therefore, to better understand the communication among neurons, we must first understand the \textquotedblleft resting\textquotedblright membrane charge. 

\begin{figure}[h!]
	\centering
	\includegraphics[width=1.0\textwidth]{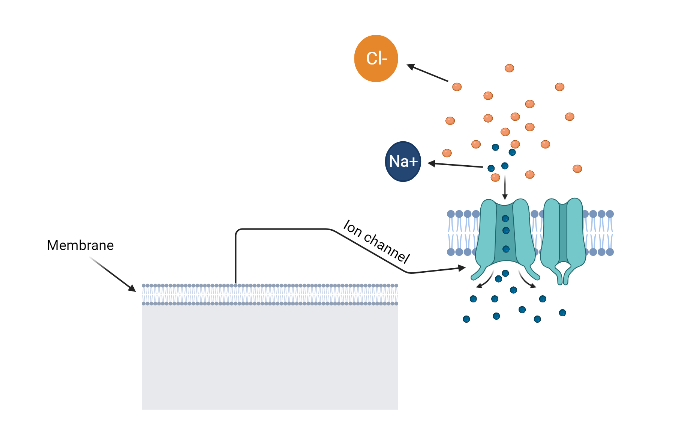}
	\caption{Schematic representation of cell's membrane potential model. At the peak action potential $\text{Cl}-$ channels close while $\text{Na}+$ channels open. $\text{Na}+$ ions leave the cell, and the membrane eventually becomes hyperpolarized (Created with BioRender.com).}
	\label{fig:membrane}
\end{figure}

\begin{figure}[h!]
	\centering
	\includegraphics[width=0.9\textwidth]{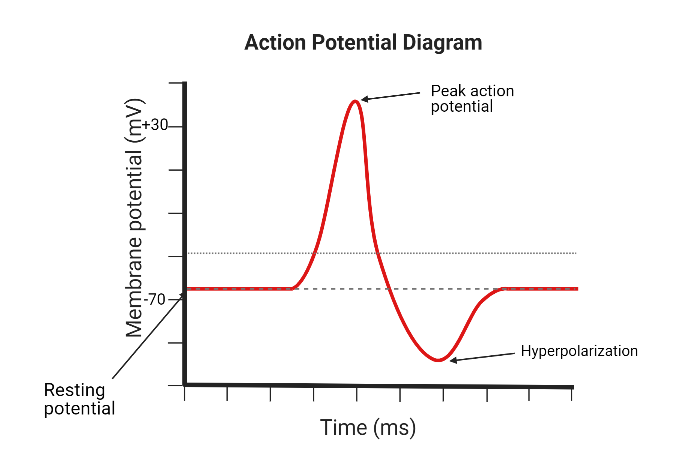}
	\caption{Schematic representation of an action potential illustrates its various phases as the action potential passes a point on a cell membrane (Created with BioRender.com). }
	\label{fig:action-potential}
\end{figure}

In general, in surrounding of each ion there is a lipid bilayer membrane that is impermeable to charged molecules or ions. Hence, ions must pass through special proteins called ion channels that span the membrane to enter or exit the neuron (see, e.g., \cite{Rye2016, Laing2010}). The activity of ion channels have different configurations: open, closed, and inactive, for instance, see in Fig. \ref{fig:membrane}. These ion channels are sensitive to the environment and can change their shape accordingly. Ion channels that change their structure in response to voltage changes are called voltage-gated ion channels. Voltage-gated ion channels regulate the relative concentrations of different ions inside and outside the cell membrane. The difference in total charge between the inside and outside of this cell is called the membrane potential. Moreover, a neuron can receive input from other neurons and the transmission of a signal between neurons is generally carried by a chemical called a neurotransmitter. Transmission of a signal within a neuron is carried by a brief reversal of the resting membrane potential called an action potential. An action potential is a brief reversal of the resting membrane potential caused by the transmission of a signal within a neuron from dendrite to axon terminal, see, e.g., Fig \ref{fig:action-potential}. In general, when neurotransmitter molecules bind to receptors located on a neuron’s dendrites, ion channels open. In fact, this opening ion channels allows positive ions to enter the neuron at excitatory synapses and results in depolarization of the membrane a decrease in the difference in voltage between the inside and outside of the neuron (see, e.g., \cite{Laing2010}). Moreover, these ion channels open only for a short time that caused more positive ions inside the membrane than outside and positively-charged ions diffuse out. As these positive ions go out, the inside of the membrane once again becomes negative with respect to the outside, the phenomenon called hyperpolarization.   

\begin{figure}[h!]
	\centering
	\includegraphics[width=1.0\textwidth]{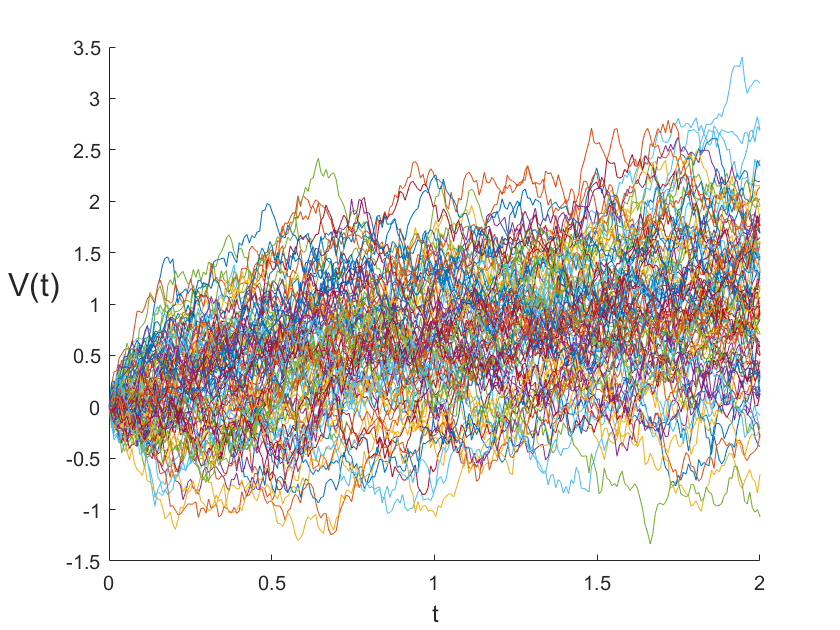}
	\caption{(Color online) Ornstein-Uhlenbeck process for $100$ sample paths with $\mu=1$ and $\sigma=1$. This figure is obtained by using Euler-Maruyama method for SDEs (see, e.g., \cite{Higham2001}).}
	\label{fig:OU}
\end{figure}

There are many different types of stochastic processes relevant to a quantitative modelling of stochastic neural activity. For instance, the behaviour of ion channels, as described above, which play an important role in the neural dynamics for action potential generation. Specifically, the membrane potential setting can vary continuously
and is driven by synaptic noise and channel noise. Hence, one of simplest models describing the membrane potential evolution is provided by the one-dimensional processes. Furthermore, the membrane fluctuations obey a simple SDE which is formally equivalent to the Ornstein-Uhlenbeck process from statistical physics, see an example in Fig \ref{fig:OU} where we plot randomly $100$ sample paths of the process using Euler-Maruyama method for SDEs (see, e.g., \cite{Higham2001}). It is clear that the drift term in the Ornstein-Uhlenbeck process represents the amount of momentum lost by drag
force and the diffusion term is subjected to a rapidly fluctuating random force. Therefore, the sample paths of this process help to explain the movement of particles in e.g. biological systems, social dynamic models, etc.

In what follows, we discuss the application of the reflected Ornstein-Uhlenbeck process in neuroscience, as well as another approach from a lattice gas setting to analyze the behavior of neurons in a membrane potential model.

Let us first recall the Ornstein-Uhlenbeck process modeling of the action potential of neurons (see, e.g., \cite{Ha2009, Xing2012, Bo2013}). 
Let $\{V_t : t \geq 0\}$ be a reflected Ornstein-Uhlenbeck process defined on $[r, \infty)$ with drift $(\mu - V_t/\gamma)$ and constant diffusion parameter $\sigma$. Then, the process $\{V_t : t \geq 0\}$ satisfies the following SDE:

\begin{align}\label{OU}
\begin{cases}
dV_t = (\mu - \frac{1}{\gamma}V_t)dt + \sigma dW_t + dL_t, \\
V_0 \in [r, \infty),
\end{cases}
\end{align}
where $W_t$ is a standard Brownian motion, while $\{L_t: t \geq 0\}$ is a continuous non-decreasing process with $L_0 = 0$. In particular, the process $\{L_t: t \geq 0\}$ increases only when $V_t = r$ and keeps $V_t \geq r$ for all $t$. In \eqref{OU}, the quantity $1/\gamma$ is the speed of reversion to its long term mean $\mu \gamma$ and $\sigma$ represents the diffusion coefficient which is instantaneous. The Ornstein-Uhlenbeck process can be interpreted as the scaling limit of mean-reverting discrete Markov chains, analogous to Brownian motion as the scaling limit of simple random walk. One of particular examples is the Ehrenfest Urn model discussed in \cite{Casas2015}. 

It is clear that after the membrane potential is reached and the neuron fires, it frequently goes below the resting baseline starting level, namely hyperpolarization. Hence, a reflecting
boundary setting for the Ornstein-Uhlenbeck model can capture the process of the membrane potential
with the hyperpolarization level. In particular, a reflecting boundary in the Ornstein-Uhlenbeck model is the maximum hyperpolarization level. By using this reflecting stochastic dynamics, we can control the action potential that goes down immediately to the resting state after reaching the peak. 

\begin{figure}[h!]
	\centering
	\includegraphics[width=0.9\textwidth]{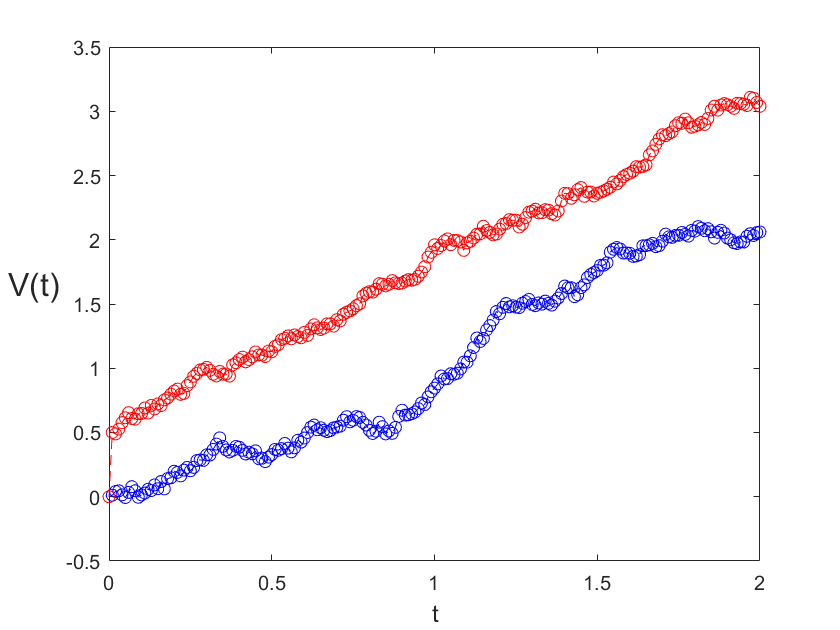}
	\caption{Membrane potential following Ornstein-Uhlenbeck process with $\mu = 1.2$, $\gamma=10$ and $\sigma=0.3$ for two cases: the red curve represents the SDEs with reflecting boundary condition at $r=0.5$, while the blue curve represents the SDEs without reflecting boundary condition.}
	\label{fig:fig2}
\end{figure}

Using the Euler-Maruyama method for SDEs (see, e.g., \cite{Higham2001}), we provide a numerical example of the membrane potential evolution in Fig. \ref{fig:fig2} following the Ornstein-Uhlenbeck process. 

We have shown a simple neuronal model based on one-dimensional processes to describe the membrane potential evolution where we are interested in the movement of only one single neuron. In general, since the human nervous system consists of billions of neural cells (or neurons), this novel ingredient brings
a lot of challenging questions to the modeling of such biological systems. For example, one of the questions we will try to answer in remainder of this paper is how to capture the behaviour of many ions when the action potential reaches peak and then goes down to the hyperpolization state.

To address this question we took the inspiration from the queueing theory model presented in Section \ref{model}. We assume that there are two different populations of ions: $\text{Na}+$ and $\text{Cl}-$ moving in a cell's membrane potential (a $2\text{D}$ lattice). We use the lattice gas dynamics, similar to Section \ref{model} for the case of $\varepsilon=0$, where now this model describes the dynamics of two populations of ions in a membrane potential. Here, $\text{Na}+$ ions can be seen as active customers, while $\text{Cl}-$ ions can be consider as passive customers. For all of the detailed analysis of reflecting stochastic dynamics of such ions ($\text{Na}+$ and $\text{Cl}-$), we use an analogous description to Section \ref{model}, replacing $A$ by $\text{Na}+$ and $P$ by $\text{Cl}-$. For the lattice gas approximation in this application, we provide the detailed setting in the next Section.

\subsection{Numerical results}

In this section, we consider a cell's membrane potential of two different populations of ions: $\text{Na}+$ and $\text{Cl}-$. The geometry is a square lattice $\Lambda:=\{1,\dots,L\}\times\{1,\dots,L\}\subset \mathbb{Z}^2$ of side 
length $L$, where $L$ is an odd positive integer number. $\Lambda$ will be referred in this context as {\em nerve net}.  
An element $x=(x_1,x_2)$ of the room $\Lambda$ is called \emph{site} 
or \emph{cell}.
Two sites $x,y\in\Lambda$ are said \emph{nearest neighbors} if and only 
if $|x-y|=1$. In addition, $N_{\text{Na}+}$ is the total number of $\text{Na}+$ ions, $N_{\text{Cl}-}$ is the total number of $\text{Cl}-$ ions with $N:=N_{\text{Na+}}+N_{\text{Cl-}}$ and $N_{\text{Na+}}, N_{\text{Cl-}}, N \in \mathbb{N}$. 

\begin{figure}[h!]
	\centering
	\begin{tabular}{lll}
		\includegraphics[width=0.33\textwidth]{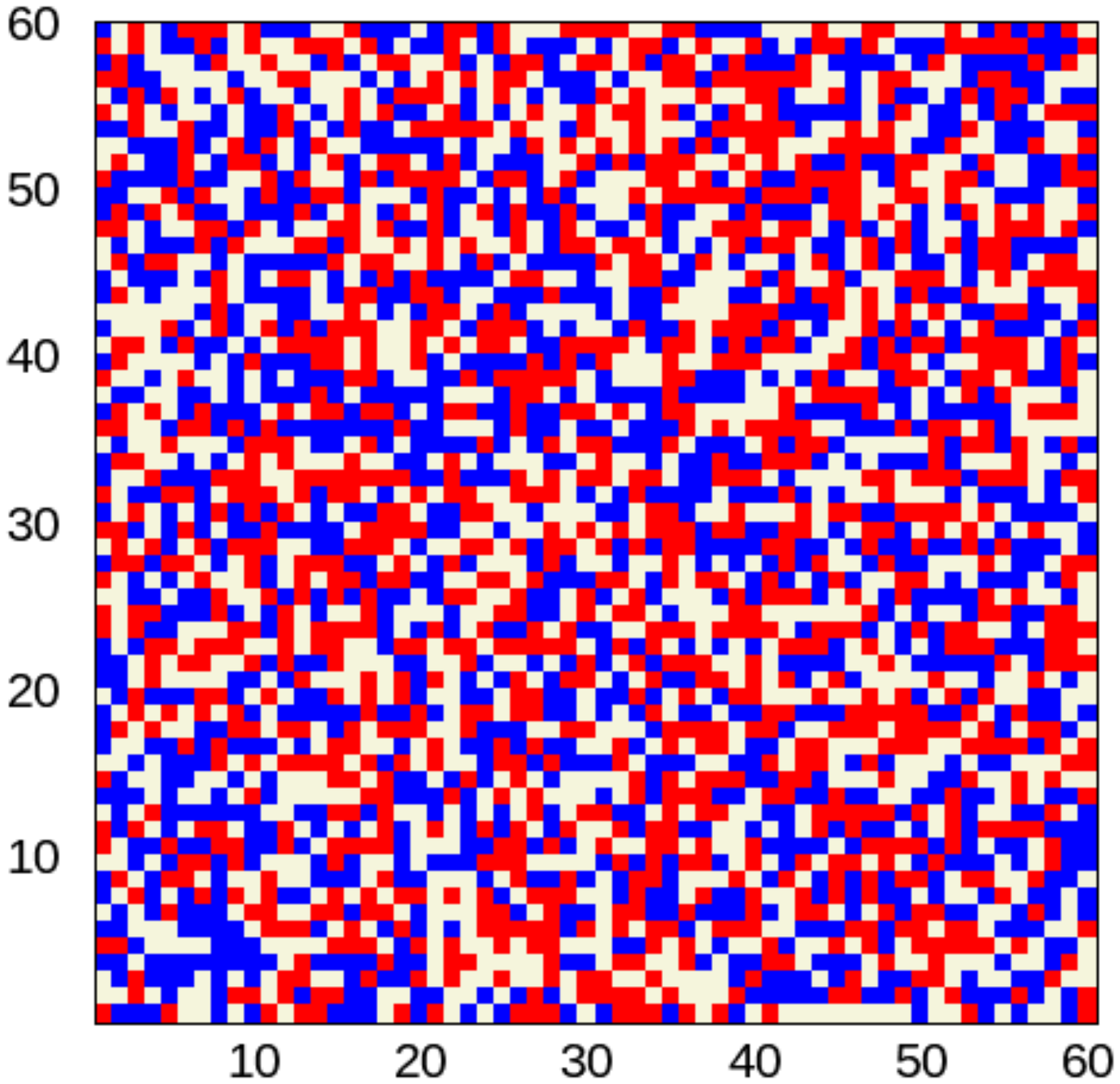}&
		\includegraphics[width=0.33\textwidth]{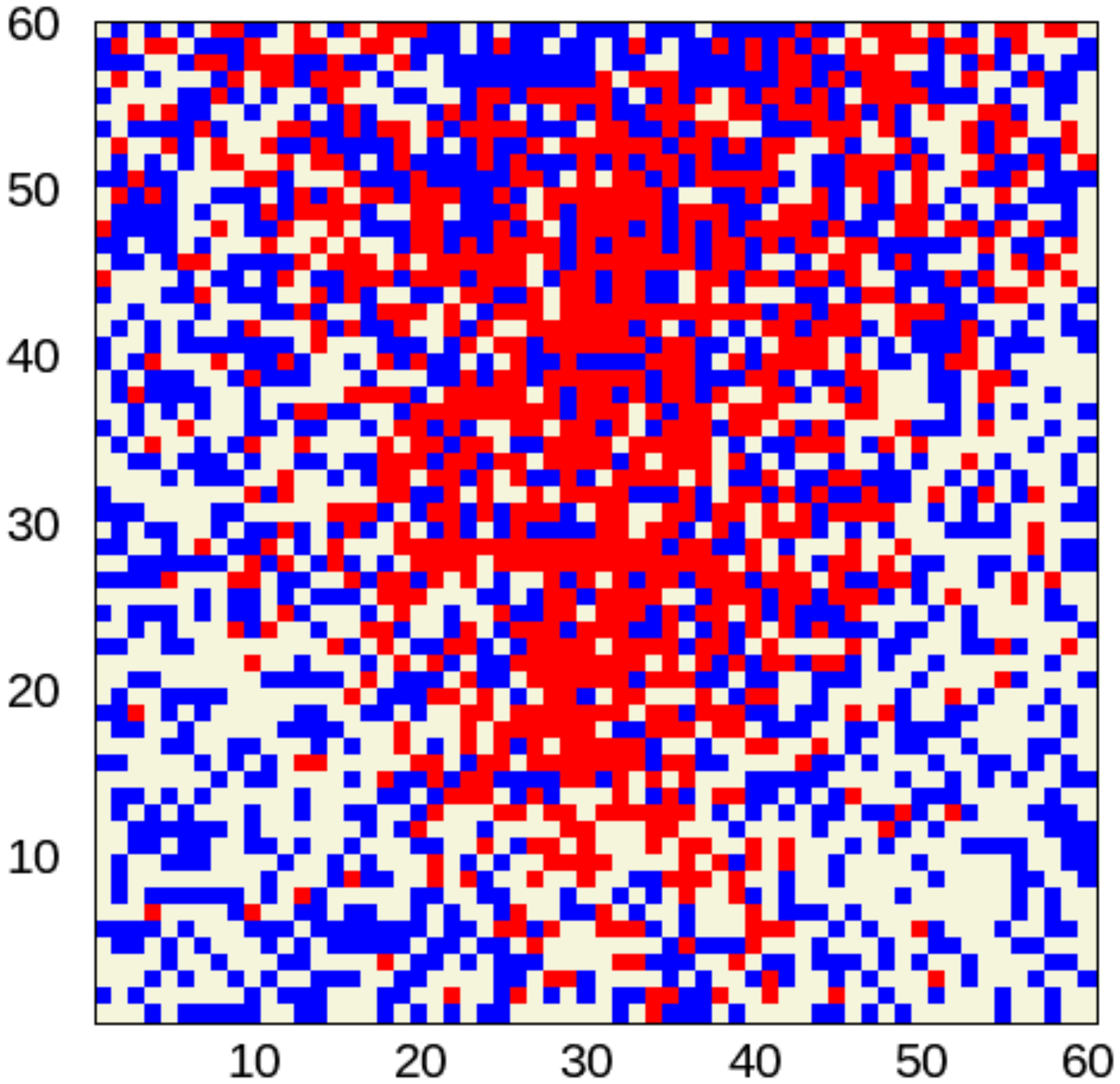}&
		\includegraphics[width=0.33\textwidth]{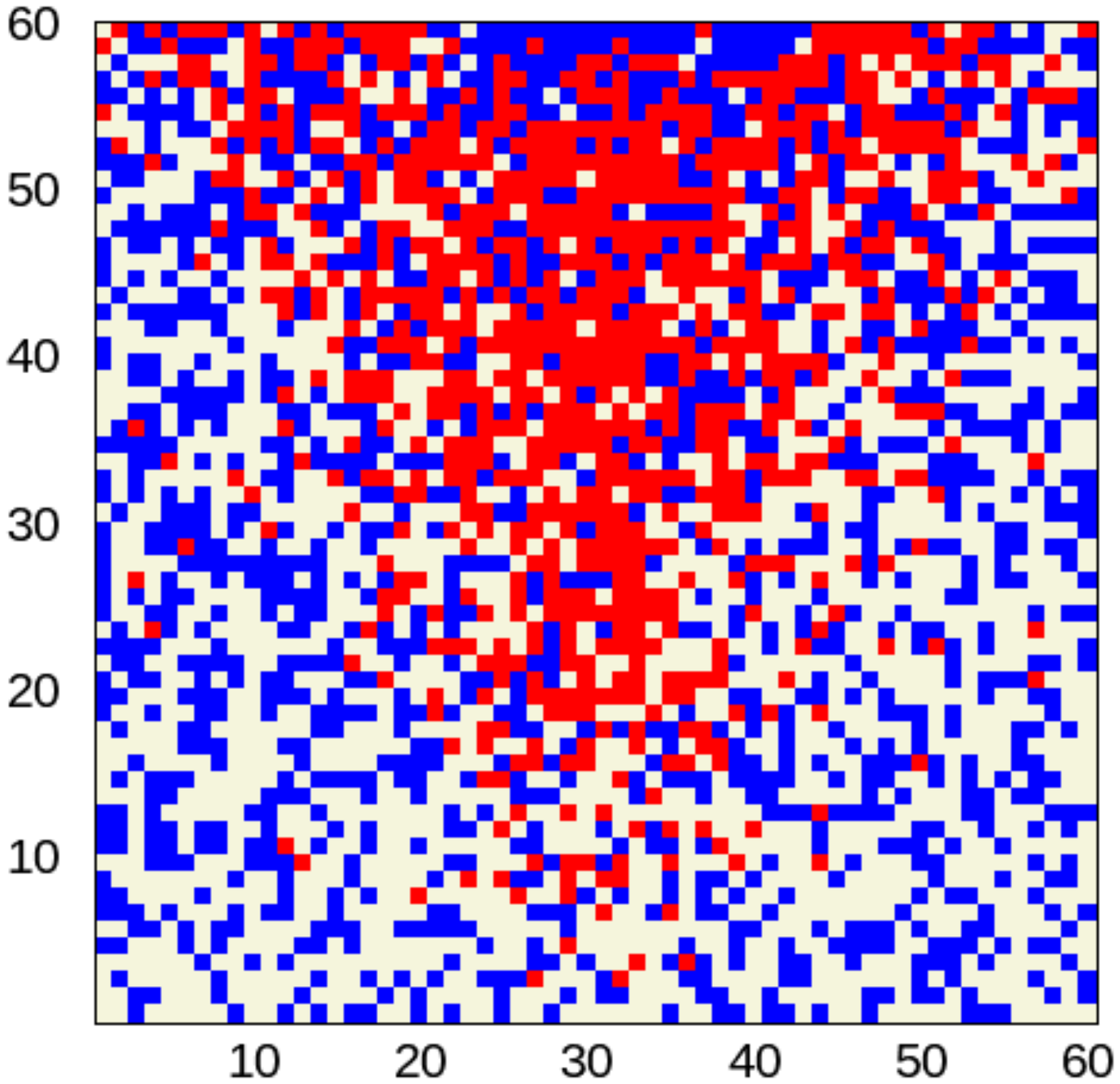}\\[0.1cm]
		\includegraphics[width=0.33\textwidth]{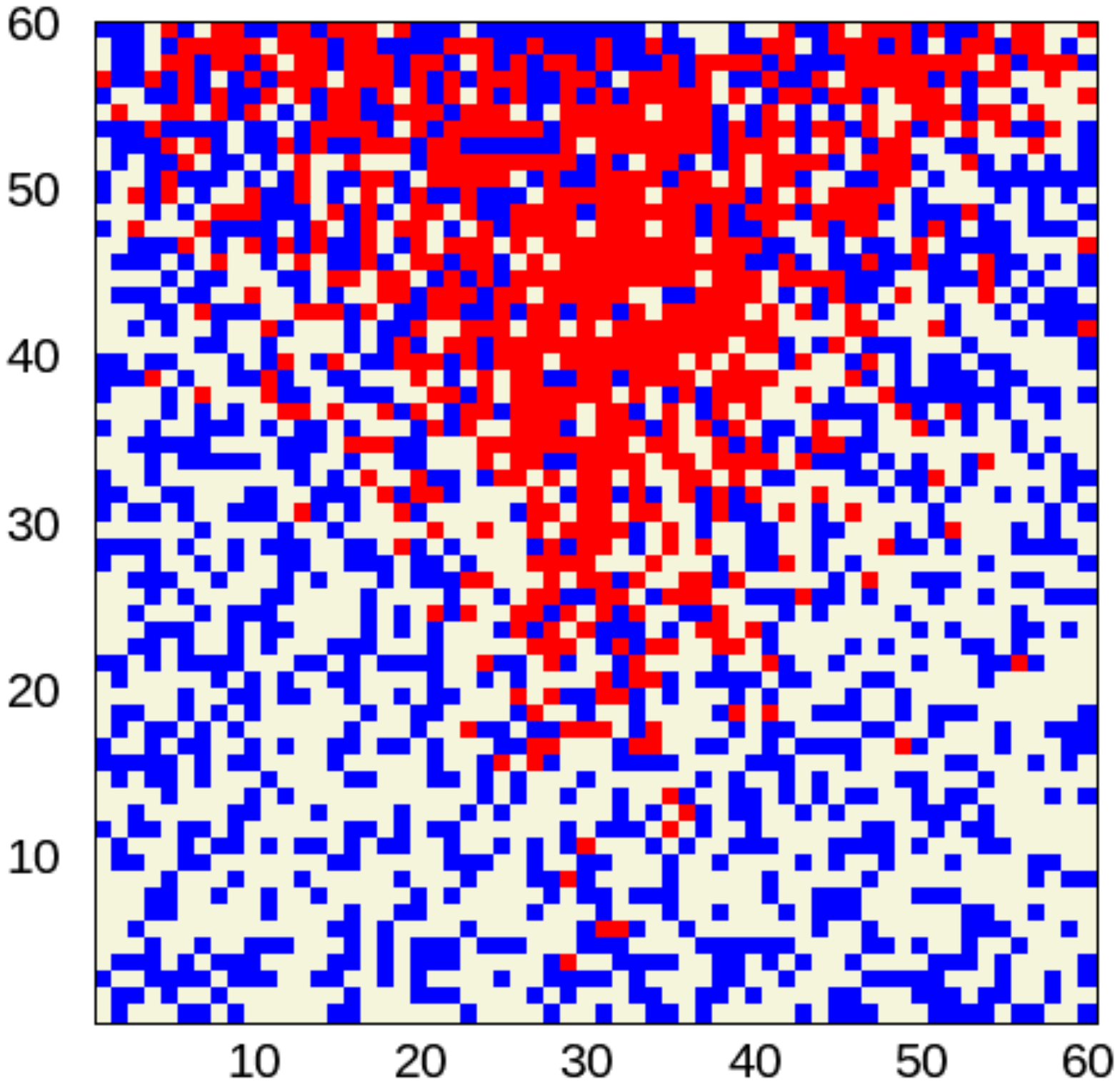}&
		\includegraphics[width=0.33\textwidth]{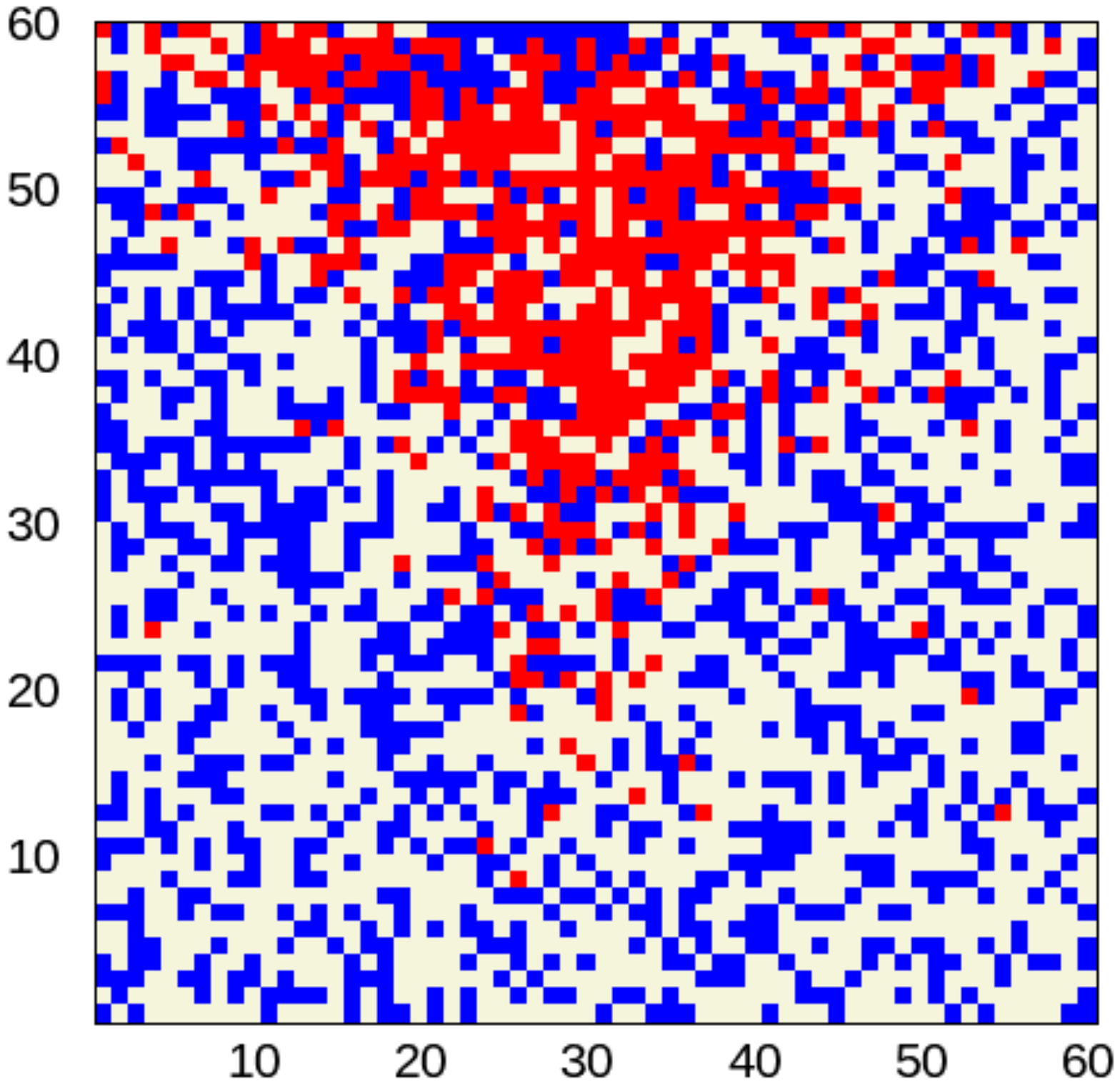}&
		\includegraphics[width=0.33\textwidth]{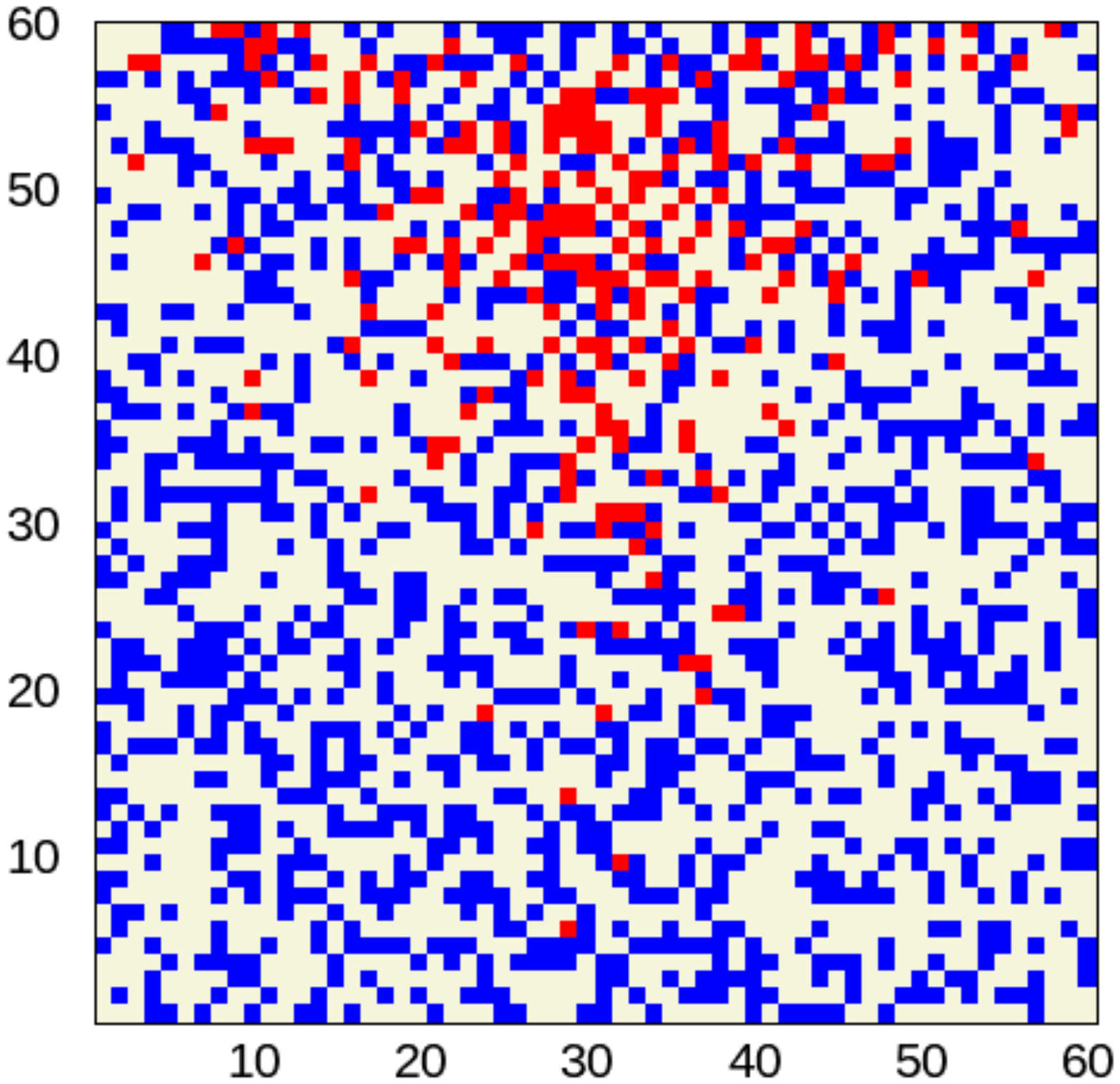}
		\\[0.1cm]
		\includegraphics[width=0.33\textwidth]{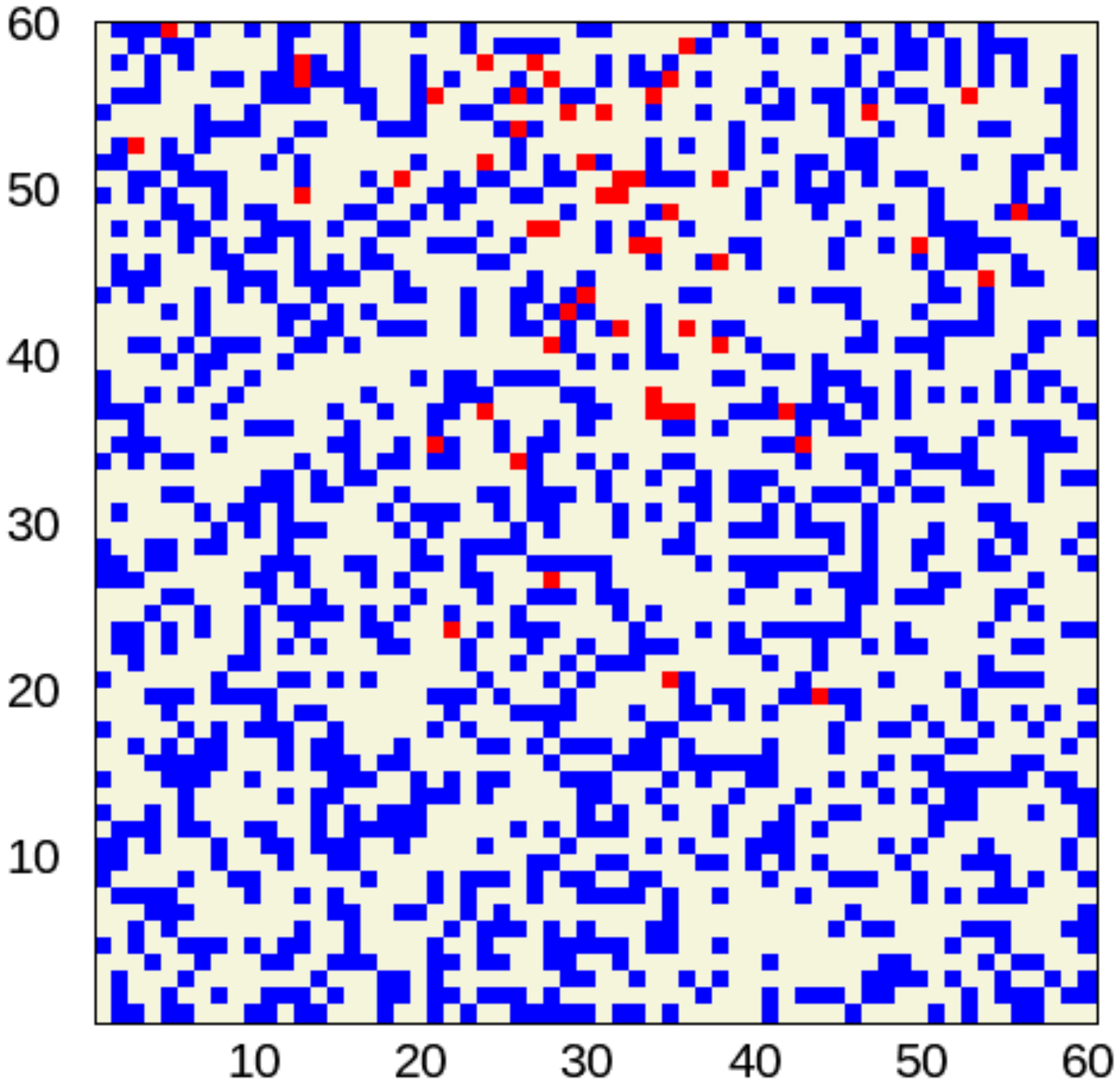}&
		\includegraphics[width=0.33\textwidth]{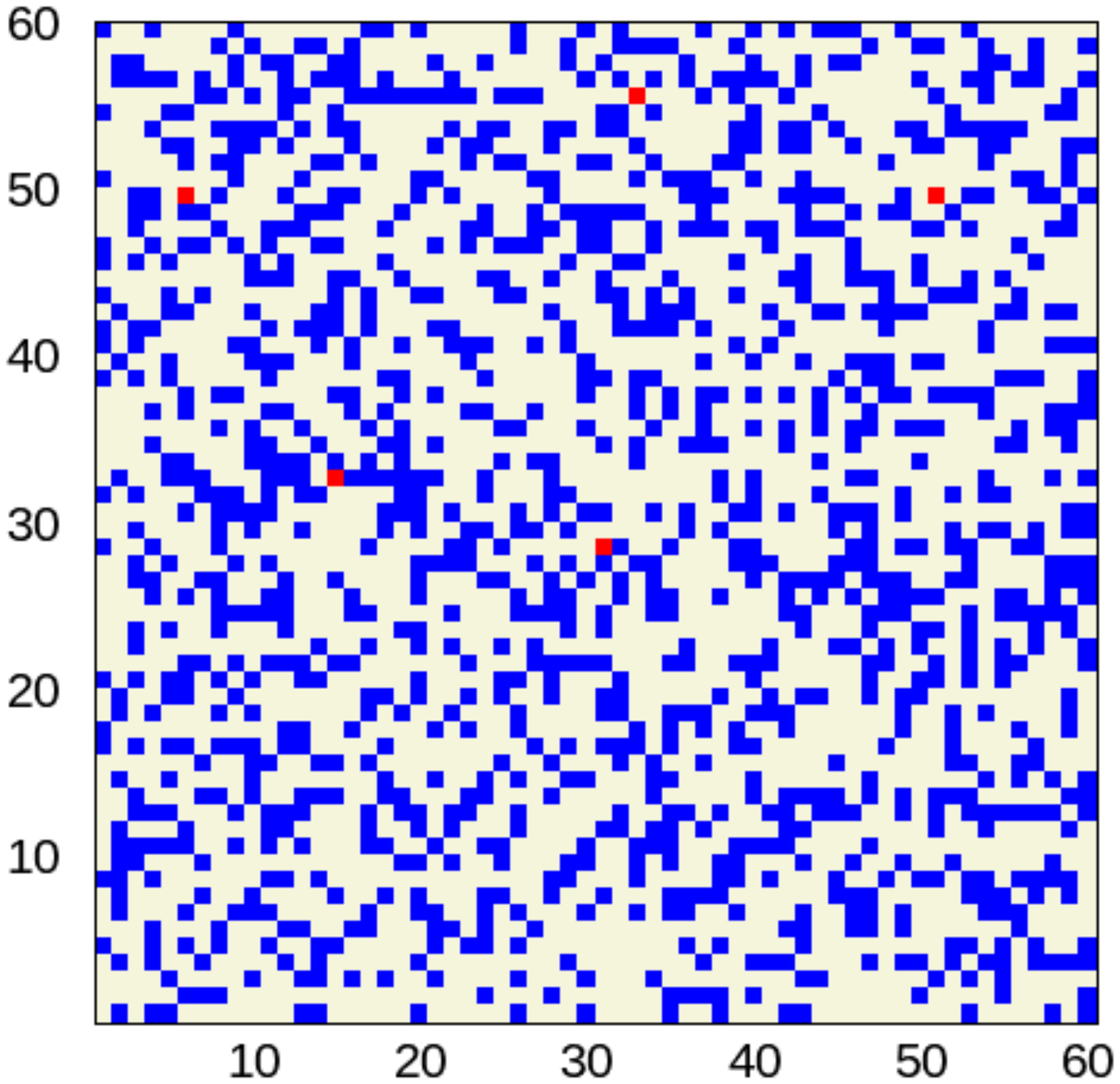}&
		\includegraphics[width=0.33\textwidth]{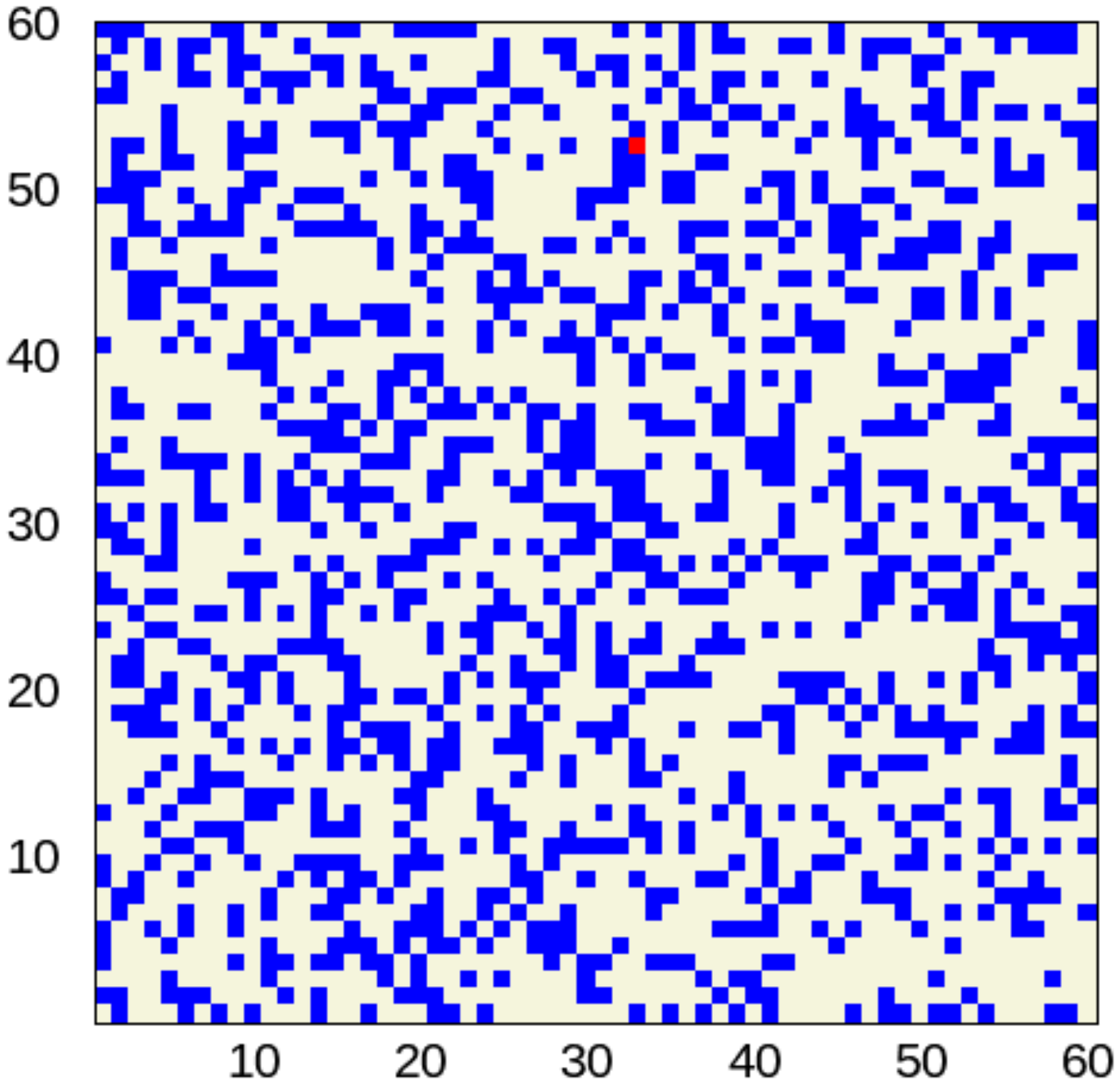}		
	\end{tabular}
	\caption{\small (Color online) Configurations of the model at different 
		times (increasing in lexicographic order).
		Parameters: $L=60$, $w_{\mathrm{ex}}=20$ and $\varepsilon=0.2$. 
		Red pixels represent active particles, blue pixels denote passive 
		particles, and gray sites are empty. 
		In the initial configuration (top left panel) there are $1200$ active and 
		$1200$ passive customers.}
	\label{fig:fig3}
\end{figure}

Furthermore, we assume that ion channels located on the top row of the geometry, and that there are $\omega$ ion channels at the exit door. The ion channels only open for $\text{Na}+$ ions to leave the cell membrane, while the $\text{Cl}-$ ions cannot leave the domain due to the inaccessibility of the channels. This model is illustrated in Fig. \ref{fig:fig3}, where we show the configuration of the system at different times. The opening of channels that let only positive ions flow out of the cell can cause the hyperpolarization. This phenomenon is observed when the membrane potential becomes more negative at a particular spot on the neuron's membrane. 
%
To simulate this process, as we have mentioned above, we exploit the idea of a lattice gas model to describe the dynamics of neurons in a cell's membrane potential model. We use the same numerical scheme as in Section \ref{KMC}, with the only difference here is that active customers are considered as $\text{Na}+$ ions and passive customers are represented by $\text{Cl}-$ ions. The neuroscience model that we study here has been developed following the idea presented in Section \ref{model} for the case of $\varepsilon=0$. This means that both $\text{Na}+$ and $\text{Cl}-$ ions perform a 
symmetric simple exclusion dynamics on the whole lattice, where only $\text{Na}+$ ions can flow out the channels. At this moment, we consider the case without chemical reactions terms in our neuronal dynamics. To get closer to a more realistic picture, the dynamics of neurons in a cell's membrane potential model would have to include chemical reactions terms. However, in such a case, the situation becomes a lot more challenging due to a feedback mechanism of the chemical reaction quantities together with interactions between the dynamics of ions and the surrounding environment. Such a generalized model is not a subject of our investigation here, and it will be reported elsewhere.

\begin{figure}[h!]
	\centering
	\begin{tabular}{ll}
		\includegraphics[width=0.45\textwidth]{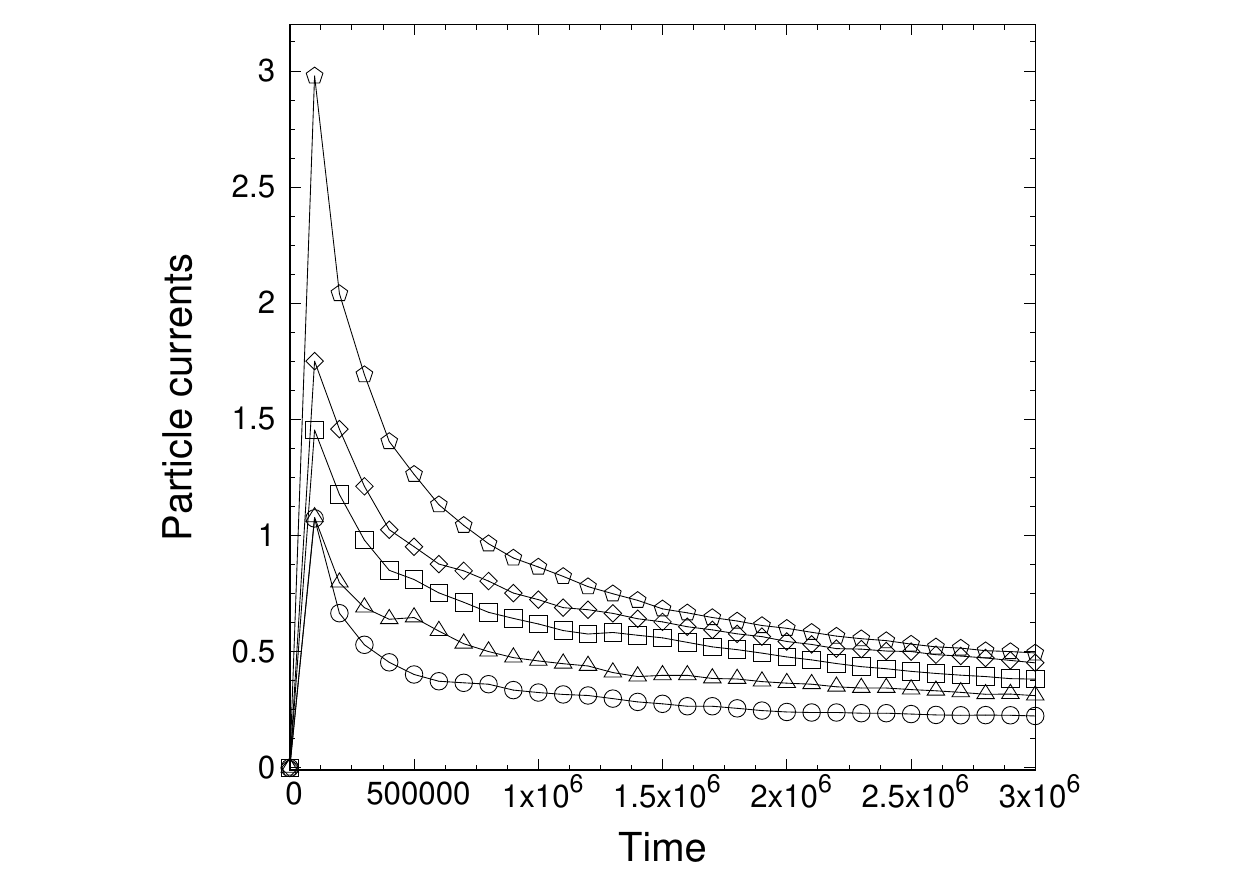}&
		\includegraphics[width=0.45\textwidth]{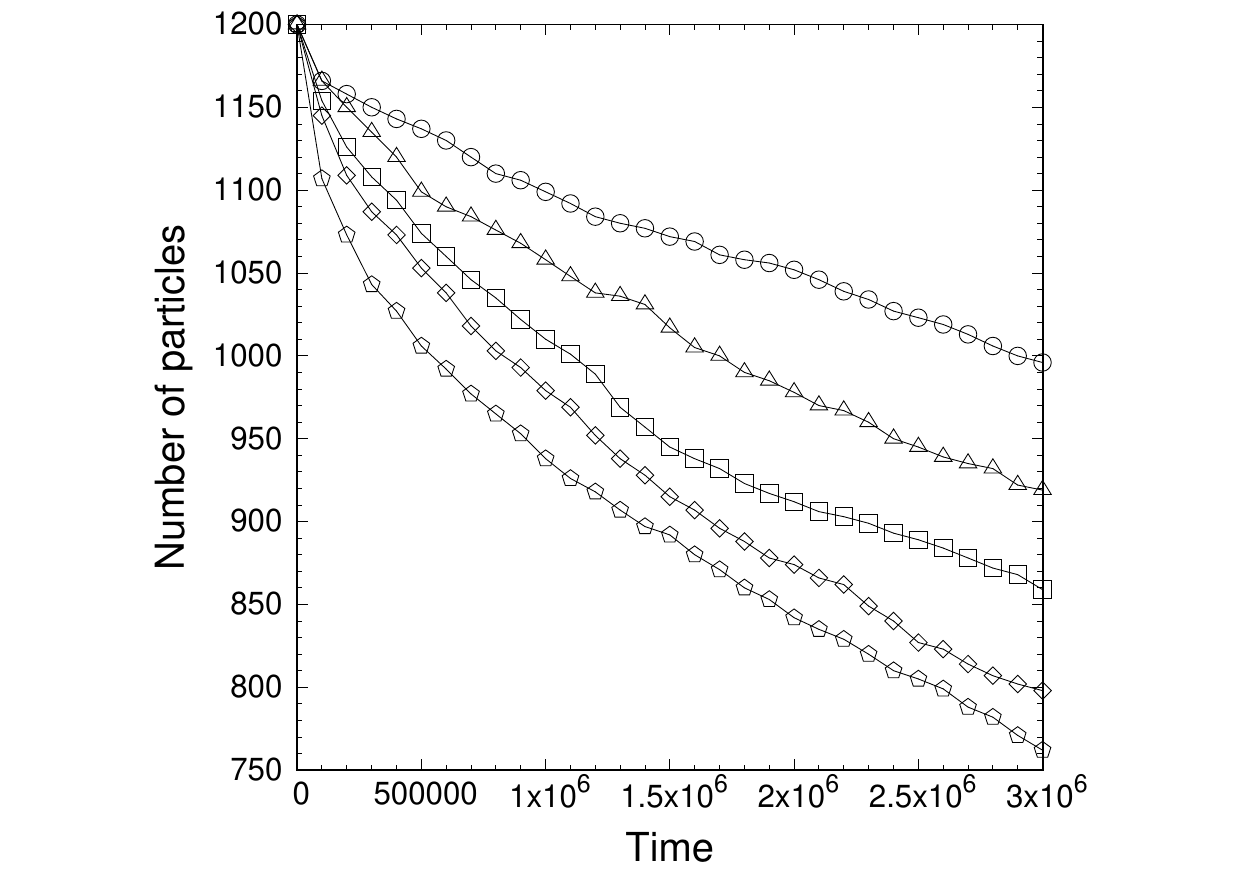}
	\end{tabular}
	\caption{Left panel: Evolution of the current as a function of time for $L=60$, $N_{\text{Na+}}=N_{\text{Cl}-}=1200$ (empty symbols), $\omega = 15$ (circles), $\omega = 20$ (triangles), $\omega = 30$ (squares), $\omega = 40$ (diamonds), $\omega = 60$ (pentagons) (circles). Right panel: Behavioral pattern of the crowd for $L=60$, $N_{\text{Na+}}=N_{\text{Cl}-}=1200$ (empty symbols), $\omega = 15$ (circles), $\omega = 20$ (triangles), $\omega = 30$ (squares), $\omega = 40$ (diamonds), $\omega = 60$ (pentagons) (circles).}
	\label{fig:fig4}
\end{figure}

\begin{figure}[h!]
	\centering
	\begin{tabular}{ll}
		\includegraphics[width=0.45\textwidth]{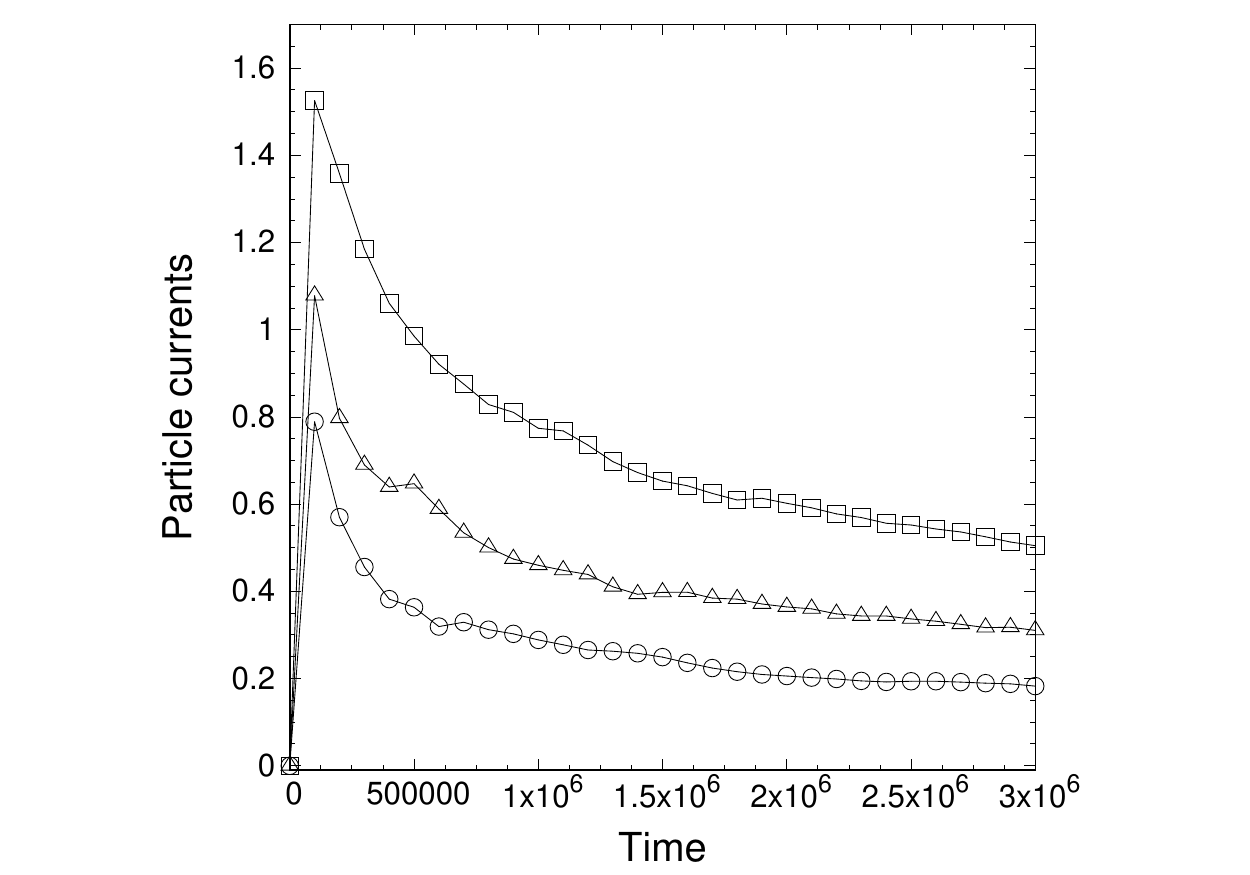}&
		\includegraphics[width=0.45\textwidth]{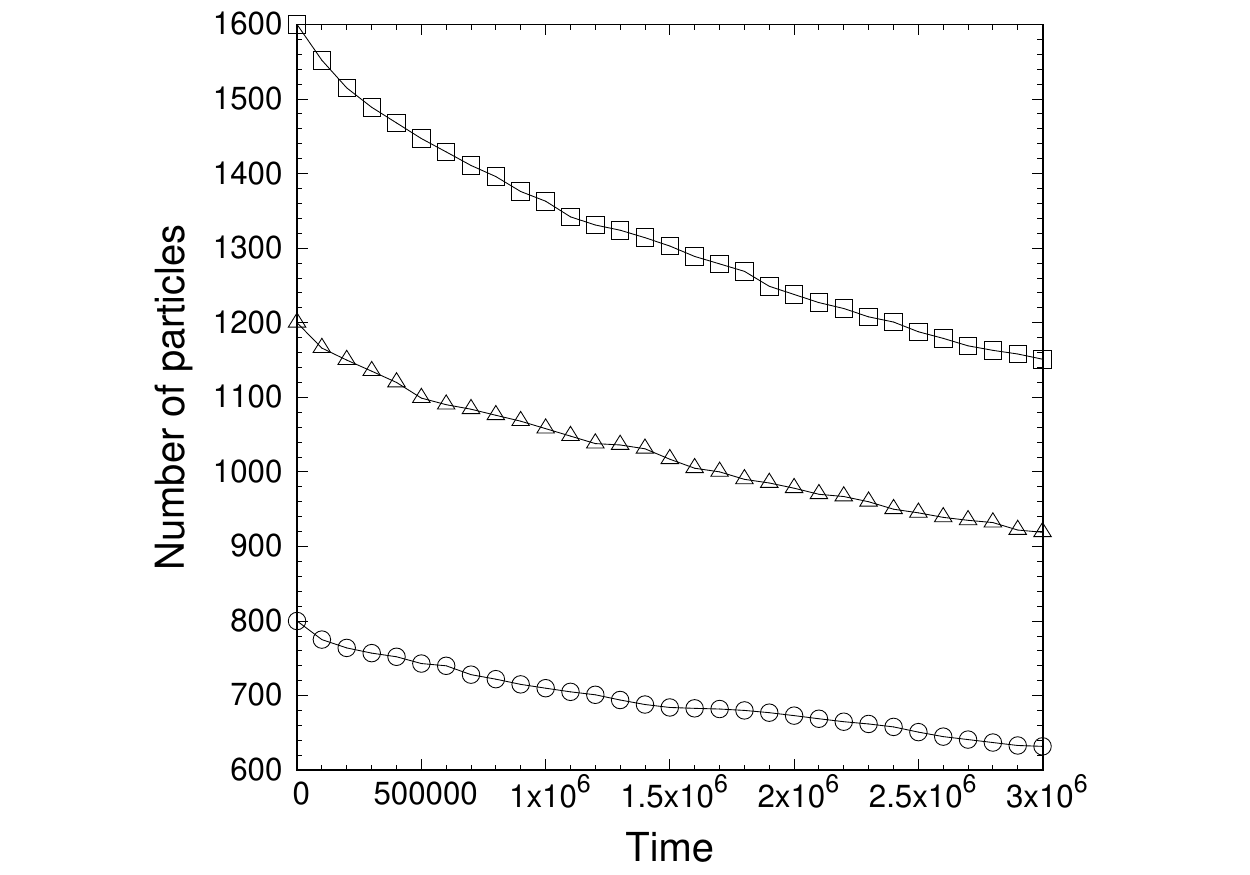}
	\end{tabular}
	\caption{Left panel: Evolution of the current as a function of time for $L=60$, $N_{\text{Na+}}=800, N_{\text{Cl}-}=1600$ (empty circles), $N_{\text{Na+}}=N_{\text{Cl}-}=1200$ (empty triangles), $N_{\text{Na+}}=1600, N_{\text{Cl}-}=800$ (empty squares). Right panel: Behavioral pattern of the crowd for $L=60$, $N_{\text{Na+}}=800, N_{\text{Cl}-}=1600$ (empty circles), $N_{\text{Na+}}=N_{\text{Cl}-}=1200$ (empty triangles), $N_{\text{Na+}}=1600, N_{\text{Cl}-}=800$ (empty squares).}
	\label{fig:fig5}
\end{figure}

The main numerical results of this part of our inverstigation is shown in Figs. \ref{fig:fig4} and \ref{fig:fig5}. We have plotted the ion currents and the behavioral pattern of the ions as a function of time, for different numbers of the ion channels in Fig. \ref{fig:fig4} and for different densities of $\text{Na}+$ and $\text{Cl}-$ ions in the cell in Fig. \ref{fig:fig5}. We observe that the current of $N_{\text{Na+}}$ ions increases when we increase the number of open ion channels in the left panel of Fig. \ref{fig:fig4}. This is visible also in the corresponding $N_{\text{Na+}}$ ion exit times in the right panel of Fig. \ref{fig:fig4}. Furthermore, we have examined also the numerical results for different densities of $N_{\text{Na+}}$ ions in the given geometry. In Fig. \ref{fig:fig4}, we obtained similar results, where we increase the density of $N_{\text{Na+}}$ ions, and the currents of ions also increase. Moreover, we have shown the behavioral pattern of $N_{\text{Na+}}$ ions at the corresponding $N_{\text{Na+}}$ ion exit times in the right panel of Fig. \ref{fig:fig5}. The activity of neurons in a cell's membrane potential problem can be predicted via this behavioral system so that we can control the currents of ions in the cell for specific purposes. Our lattice model opens an interesting research direction where the behavior of many ions can be investigated in the study of neuroscience by developing further its intrinsic links with reflecting stochastic processes of mixed population dynamics. The next step would be the investigation of the dynamics of ions in the presence of a feedback mechanism of chemical reaction terms in the associated neuronal models.
\section{Concluding remarks}\label{concluding-remarks}

We have analyzed numerically the reflecting stochastic dynamics of active-passive populations and provided two representative examples in a queueing theory model and neuroscience. We have used a statistical-mechanics-based lattice gas framework where we employ a kinetic
Monte Carlo procedure for the implementation of our corresponding models. In the first model, based on our numerical experiments, we have observed the impact of passive customers on the residence times of the active population, which allowed us to conclude that the presence of passive customers in the system increases the waiting time of the active customers. In reality, setting a limit on the presence of passive customers would allow for an optimization approach on the waiting time of the queues. In the current consideration, we have limited our representative examples to a situation where the interaction between active and passive customers is subject to simple exclusion rules. It would be instructive to analyze the situations where nonlocal interactions among the crowd participants are also allowed. In the second model, based on our numerical experiments, we have examined the behavior of many ions when the action potential reaches peak and then go down to the hyperpolization state. Via this behavioral system, the activity of neurons in a cell's membrane potential problem can be predicted. In the current model, we investigated the behavior of many ions in the absence of chemical reaction terms. It would be interesting to extend the model for a more complex scenario where chemical reaction terms are added to the system. A comparison between the theoretical prediction with the real data would benefit further progress in the development of the presented framework.

The class of reflecting stochastic dynamics in operation research and neuroscience models provides a bridge between reflecting stochastic dynamics in a confined domain and a lattice gas approximation via a limit theorem. The ideas presented in this contribution may be extended to queueing theory models in healthcare systems subjected to an epidemic, e.g. by including symptomatic-asymptomatic populations into account \cite{Meares2020,Tadic2020,Tadic2021}, as well as to other application areas where we have to deal with uncertainties intrinsic to two groups of populations with distinct behavioral patterns.  
\section*{Acknowledgment}
Authors are grateful to the NSERC and the CRC Program for their
support. RM is also acknowledging support of the BERC 2018-2021 program and Spanish Ministry of Science, Innovation and Universities through the Agencia Estatal de Investigacion (AEI) BCAM Severo Ochoa excellence accreditation SEV-2017-0718 and the Basque Government fund AI in BCAM EXP. 2019/00432.

\bibliography{mybibfile}

\end{document}